\documentclass[12pt]{article}

\usepackage{epsfig}
\usepackage{amsfonts}
\usepackage{amssymb}
\usepackage{amsmath}
\usepackage{amsthm}
\usepackage{latexsym}
\usepackage{color}

\newtheorem{theorem}{Theorem}[section]
\newtheorem{lemma}[theorem]{Lemma}

\newtheorem{assumption}{Assumption}

\DeclareSymbolFont{AMSb}{U}{msb}{m}{n}
\DeclareMathSymbol{\N}{\mathbin}{AMSb}{"4E}
\DeclareMathSymbol{\Z}{\mathbin}{AMSb}{"5A}
\DeclareMathSymbol{\R}{\mathbin}{AMSb}{"52}
\DeclareMathSymbol{\Q}{\mathbin}{AMSb}{"51}
\DeclareMathSymbol{\I}{\mathbin}{AMSb}{"49}
\DeclareMathSymbol{\C}{\mathbin}{AMSb}{"43}

\title{Coarse-grained modeling of multiscale diffusions: the p-variation estimates.}

\author{Anastasia Papavasiliou\footnote{Department of Statistics, University of Warwick, Coventry, CV4 7AL, UK.}}

\begin{document}
\maketitle

\begin{abstract}
We study the problem of estimating parameters of the limiting equation of a multiscale diffusion in the case of averaging and homogenization, given data from the corresponding multiscale system. First, we review some recent results that make use of the maximum likelihood of the limiting equation. In particular, it has been shown that in the averaging case, the MLE will be asymptotically consistent in the limit while in the homogenization case, the MLE will be asymptotically consistent only if we subsample the data. Then, we focus on the problem of estimating the diffusion coefficient. We suggest a novel approach that makes use of the total $p$-variation, as defined in \cite{Terry book} and avoids the subsampling step. The method is applied to a multiscale OU process.
\end{abstract}


\bigskip
{\bf AMS subject classifications:}

\section{Introduction}
It is often the case that the most accurate models for physical systems are large in dimension and multiscale in nature. One of the main tasks for applied mathematicians is to find coarse-grained models of smaller dimension that can effectively describe the dynamics of the system and are efficient to use (see, for example \cite{MTV 1, MTV 2, KMS 1, KMS 2}). Once such a model is chosen, its free parameters are estimated by fitting the model to the existing data. Here, we study the challenges of this statistical estimation problem, in particular for the case where the coarse-grained model is a diffusion. Apart from the usual challenges of parameter estimation for diffusions, an additional problem that needs to be addressed in this setting is that of the mismatch between the full multiscale model that generated the data and the coarse-grained model that is fitted to the data. A first discussion of this issue, in the context of averaging and homogenization for multiscale diffusions, can be found in \cite{PS, PPS, ABT} .

A similar statistical estimation problem arises in the context of ``equation-free'' modeling. In this case, coarse-grained equations exist only locally and are locally fitted to the data. The main idea of ``equation-free'' modeling is to use these locally fitted coarse-grained equations in combination with a global algorithm (for example, Newton-Raphson) in order to answer questions about the global dynamics of the coarse-grained model (for example, finding the roots of the drift). In this process, we go through the following steps: we simulate {\it short} paths of the system for given initial conditions. These are used to locally estimate the effective dynamics. Then, we carefully choose the initial conditions for the following simulations so that we reach an answer to whatever question we set on the global dynamics of the system, as quickly and efficiently as possible (see \cite{Kevrekidis}). The statistical inference problem is similar to the one before: we have the data coming from the full model, we have a model for the effective local dynamics and we want to fit the data to this model. However, there is also an important difference: the available data is {\it short} paths of the full model. This issue has not been addressed in \cite{PS, PPS} or \cite{ABT}, where it is assumed that the time horizon is either fixed or goes to infinity at a certain rate. We will address this problem in section \ref{my stuff}, by letting the time horizon $T$ be of order ${\mathcal O}\left( \epsilon^\alpha \right)$, where $\epsilon$ is the scale separation variable and $\alpha>0$. Another important issue that we will address here is that of estimating the scale separation variable $\epsilon$.

We will focus on a very simple Ornstein-Uhlenbeck model whose effective dynamics can be described by a scaled Brownian motion. This will allow us to perform precise computations, reach definite conclusions and build our intuition about the behavior of more general diffusions. We will only tackle the homogenization case and our goal will be to estimate the diffusion coefficient of the effective dynamics. This problem has also been addressed in \cite{PS, ABT}. In both these papers, the diffusion coefficient is constant. In fact, in \cite{ABT} the authors also focus on the Ornstein-Uhlenbeck model. Our main contribution is to demonstrate that in order to compute the diffusion coefficient, one should not use the quadratic variation commonly defined as a limit where we let the size of a partition go to zero but rather as a supremum over all partitions. This definition is discussed in \cite{Terry book} and is at the core of the theory of rough paths, as it gives rise to a topology with respect to which the It\^o map is continuous. 

In section \ref{review}, we review some of the core results for multiscale diffusions and their coarse-grained models. Then, we will review the results of \cite{PS, PPS} and \cite{ABT}. Finally, We will give a more precise description of ``equation-free'' modeling.

In section \ref{my stuff}, we go on to define a new set of estimators for the diffusion parameter of the coarse-grained model, in the case of homogenization. We perform explicit computations of their $L_2$-error, which allows us to attest their performance. We conclude that they outperform the subsampled quadratic variance estimate studied in \cite{PS, ABT}. Finally, we describe a heuristic way of estimating the scale separation parameter $\epsilon$. 

\section{MLE for multiscale diffusions: A review}
\label{review}

In this section, we review some of the main concepts that come into play in multiscale modeling. First, we describe the limiting equations for multiscale stochastic differential equations. These allow us to reduce the dimension of the model. Then, we discuss the problem of the statistical estimation of parameters of the limiting equation given multiscale data and how this mismatch between model and data affects the result. Finally, we discuss a numerical algorithm that is applied when the limiting equations are completely unknown, which comes under the name of ``equation-free'' modeling.

\subsection{Limiting equations for multiscale diffusions}
\label{limiting equations}

Multiscale diffusions are a combination of two basic types of multiscale stochastic differential equations. The first is described by the following equations
\begin{equation}
\label{eq:av}
\begin{array}{ccc}
dX_t &=& f_1(X_t,Y_t)dt + \sigma_1(X_t,Y_t) dW_t \\ & & \\
dY_t &=& \frac{1}{\epsilon^2}f_2(X_t,Y_t)dt + \frac{1}{\epsilon}\sigma_2(X_t,Y_t) dV_t
\end{array}
\end{equation}
where $X_t \in {\mathcal X}$ and $Y_t \in {\mathcal Y}$ and ${\mathcal X}, {\mathcal Y}$ are finite dimensional Banach spaces. We call $X$ the slow variable, $Y$ the fast variable and $\epsilon$ the scale separation parameter. The main assumptions are the following:
\begin{assumption}
\label{ass:av}
\begin{itemize}
\item[(i)] The solution of the system exists.
\item[(ii)] The equation 
\[ dY^x_t =\frac{1}{\epsilon^2}f_2(x,Y^x_t)dt + \frac{1}{\epsilon}\sigma_2(x,Y^x_t) dV_t \]
is ergodic with unique invariant measure $\mu_x$, for every $x\in{\mathcal X}$.
\end{itemize}
\end{assumption}
We expect that by the time $X$ takes a small step $\Delta \sim{\mathcal O}\left( 1 \right)$, 
\[ \frac{1}{\Delta}\int_t^{t+\Delta} f_1(X_s,Y_s)ds \approx \int f_1(X_t, y) \mu_{X_t}(dy) \]
as a result of the ergodicity of $Y$. Similarly,
\[ \frac{1}{\Delta}\int_t^{t+\Delta} \sigma_1(X_s,Y_s)\sigma_1(X_s,Y_s)' ds \approx \int \sigma_1(X_t, y)\sigma_1(X_t, y)' \mu_{X_t}(dy). \]
where by $(\cdot)'$ we denote the transpose of a vector. We set
\[ \bar{f}_1 (x) = \int f_1(x, y) \mu_{x}(dy),\ \ \bar{\sigma}_1 (x) = \left( \int \sigma_1(x, y)\sigma_1(x, y)' \mu_{x}(dy) \right)^\frac{1}{2}\]
and
\begin{equation}
\label{averaged eq} 
d\bar{X}_t = \bar{f}_1(\bar{X}_t)dt + \bar{\sigma}_1(\bar{X}_t) dW_t.
\end{equation}
We call (\ref{averaged eq}) the averaged limiting equation and we call $\bar{X}$ the averaged limit. We expect that $X_t \approx \bar{X}_t$, provided that they have the same initial conditions. Indeed, the following holds
\begin{theorem}[\cite{PPS}]
\label{th:av} Let ${\mathcal X} = {\mathbb T}^\ell$ and ${\mathcal Y} = {\mathbb T}^{d-\ell}$. We assume that all coefficients in (\ref{eq:av}) are smooth in both $x$ and $y$ and that the matrix $\Sigma_2(x,y) = \sigma_2(x,y)\sigma_2(x,y)'$ is positive definite, uniformly in $x$ and $y$. Also, there exists a constant $C>0$ such that
\[ \left\langle z, B(x,y) z \right\rangle \geq C |z|^2,\ \ \forall (x,y)\in{\mathcal X}\times{\mathcal Y}\ \ {\rm and}\ \ z\in {\mathbb R}^{d-\ell},\]
where $\left\langle \cdot, \cdot \right\rangle$ denotes the Euclidean inner product. Then, if $X_0 = \bar{X}_0$,
\[ X \Rightarrow \bar{X}\ \ {\rm in}\ \ {\mathcal C}\left( [0,T], {\mathcal X}\right).\]
\end{theorem}
Different types of convergence have also been proven under different assumptions (see \cite{PS book} and \cite{FW book}). 

The second basic type of multiscale stochastic differential equation is described by the following equations
\begin{equation}
\label{eq:hom}
\begin{array}{ccc}
dX_t &=& \frac{1}{\epsilon}f_1(X_t,Y_t)dt  \\ & & \\
dY_t &=& \frac{1}{\epsilon^2}f_2(X_t,Y_t)dt + \frac{1}{\epsilon}\sigma_2(X_t,Y_t) dV_t
\end{array}
\end{equation}
where $X_t \in {\mathcal X}$ and $Y_t \in {\mathcal Y}$ and ${\mathcal X}, {\mathcal Y}$ are finite dimensional Banach spaces. As before, we call $X$ the slow variable, $Y$ the fast variable. In addition to assumption \ref{ass:av}, we assume that
\begin{assumption}
\label{ass:hom}
\[ \int_{\mathcal Y} f_1 (x,y) \mu_x(dy) = 0,\ \ \forall x\in{\mathcal X}\]
where $\mu_x$ as defined in assumption \ref{ass:av}.
\end{assumption}
Then, we expect that by the time $X$ takes a small step $\Delta \sim{\mathcal O}\left( 1 \right)$, 
\begin{eqnarray*} 
\frac{1}{\Delta \epsilon}\int_t^{t+\Delta} f_1(X_s,Y_s)ds \approx \frac{1}{\Delta \epsilon}\int_t^{t+\Delta} f_1(X_t,Y^{X_t}_s)ds  
\end{eqnarray*}
It follows from the Central Limit Theorem for ergodic Markov Processes (see \cite{Chen book}) that this will converge to a random number. More precisely, let us set
\[ \bar{f}_1 (x) = \int_{\mathcal X} \int_0^\infty f_1(x,y) \left(P_s \partial_x f_1(x,\cdot)\right)(y)'\mu_x(dy), \]
and
\[ \bar{\tau} (x) = \left( 2 \int_{\mathcal X} \int_0^\infty f_1(x,y) \left(P_s f_1(x,\cdot)\right)(y)'\mu_x(dy)\right)^\frac{1}{2}, \]
where $P_t$ are the transition kernels of the diffusion $Y^x$. Finally, we set
\begin{equation}
\label{homogenized eq}
d\bar{X}_t = \bar{f}_1(\bar{X}_t)dt + \bar{\tau}(\bar{X}_t) dW_t.
\end{equation}
We call $\bar{X}$ the homogenized limiting equation. As before, we expect that $X_t \approx \bar{X}_t$, provided that they have the same initial conditions. Indeed, similar to the averaging case, we can prove the following:
\begin{theorem}[\cite{PPS}]
\label{th:hom} Let ${\mathcal X} = {\mathbb T}^\ell$ and ${\mathcal Y} = {\mathbb T}^{d-\ell}$. We assume that all coefficients in (\ref{eq:hom}) are smooth in both $x$ and $y$ and that the matrix $\Sigma_2(x,y) = \sigma_2(x,y)\sigma_2(x,y)'$ is positive definite, uniformly in $x$ and $y$. Also, there exists a constant $C>0$ such that
\[ \left\langle z, B(x,y) z \right\rangle \geq C |z|^2,\ \ \forall (x,y)\in{\mathcal X}\times{\mathcal Y}\ \ {\rm and} z\in {\mathbb R}^{d-\ell},\]
where $\left\langle \cdot, \cdot \right\rangle$ denotes the Euclidean inner product. Then, if $X_0 = \bar{X}_0$ and assumption \ref{ass:hom} holds, we get that
\[ X \Rightarrow \bar{X}\ \ {\rm in}\ \ {\mathcal C}\left( [0,T], {\mathcal X}\right).\]
\end{theorem}
Again, different types of convergence have also been proven under different assumptions (see \cite{PS book} and \cite{FW book}). 

Theorems \ref{th:av} and \ref{th:hom} allow us to replace the $(X_t,Y_t)$ system by $\bar{X}_t$. If we are only interested in the slow dynamics of the system, this allows us to reduce the dimension of the problem. For example, using the limiting equations we can simulate the slow dynamics of the process much faster, not only because of the dimension reduction but also because the dynamics of $\bar{X}$ do not depend on $\epsilon$. Thus, the step of any numerical algorithm used to simulate the dynamics can be of order ${\mathcal O}\left(1 \right)$ rather than ${\mathcal O}\left(\epsilon^2 \right)$ which would have been the case if we wanted to simulate the full multiscale system.

\subsection{Parameter estimation for multiscale diffusions: a review}
\label{PE review}

The theory reviewed in section \ref{limiting equations} allows us to reduce the dimension of a multiscale system, approximating the slow dynamics by an diffusion of smaller dimension that does not have a multiscale structure anymore. In addition to multiscale diffusions, similar results hold for ordinary and partial differential equations (see \cite{PS book}). 

It is often the case that the dynamics of the full multiscale system -- and consequently those of the limiting system -- are not completely known. For example, in the case of multiscale diffusions, the drift and variance of the full system and thus the limiting system might depend on unknown parameters. This poses a statistical problem: how can we estimate these parameters give the multiscale data? In fact, it is even more realistic to ask to find the drift and diffusion coefficient of $\bar{X}$ given only $X$. This problem has been discussed in \cite{PS,PPS,ABT}.

More precisely, in \cite{PS}, the authors discuss the case where the drift of the limiting equation depends linearly on the unknown parameter while the diffusion parameter is constant. In \cite{PPS}, the authors extended the results of \cite{PS} for generic drift but did not discuss the problem of estimating the diffusion parameter. Finally, in \cite{ABT}, the authors extend the results in \cite{PS} by also proving the asymptotic normality of the estimators, but the limit their study to the Ornstein-Uhlenbeck system.
The approach taken so far is the following:
\begin{itemize}
\item[(i)] We pretend that the data comes from the limiting equation and we write down the corresponding maximum likelihood estimate (MLE) for the unknown parameters;
\item[(ii)] we study whether the mismatch between model and data leads to errors and, if so, we try to find a way to correct them. It has been shown that in the limit as the scale separation parameter $\epsilon\rightarrow 0$, the MLE corresponding to the averaged equation is consistent. However, this is not true in the case of homogenization. The method used so far to correct this problem has been that of subsampling the data by a parameter $\delta$. Then, for $\delta\sim{\mathcal O}\left( \epsilon^\alpha \right)$ and $\alpha\in [0,2]$, it has been shown that the MLE that corresponds to the homogenized equation will be consistent in the limit $\epsilon\rightarrow 0$. Also, an effort has been made to identify the optimal subsampling rate, i.e. the optimal $\alpha$. However, since $\epsilon$ is usually an unknown, this is of little practical value.
\end{itemize}
Note that a separate issue is that of writing the maximum likelihood of the limiting diffusion, which in the general multi-dimensional case can still be challenging (see \cite{Bishwal, Kutoyants}). We will not discuss this issue here, however. 

We summarize the main results for the parameter estimation of the limiting equations of multiscale diffusions in the following theorems:

\begin{theorem}[Drift estimation, averaging problem]
\label{thm: drift av}
Suppose that $\bar{f}_1$ in (\ref{averaged eq}) depends on unknown parameters $\theta$, i.e. $\bar{f}_1(x) = \bar{f}_1(x;\theta)$. Let $\hat{\theta}(x;T)$ be the MLE of $\theta$ corresponding to equation (\ref{averaged eq}). Suppose that we observe $\{ X_t, t\in[0,T] \}$ of system (\ref{eq:av}) corresponding to $\theta = \theta_0$. Then, under appropriate assumptions described in \cite{PPS} (Theorem 3.11), it is possible to show that
\[ \lim_{\epsilon\rightarrow 0} dist \left( \hat{\theta}(X;T), \theta_\epsilon \right)= 0,\ \ {\rm in\ probability} \]
where $dist\left(\cdot,\cdot \right)$ is the asymmetric Hausdorff semi-distance and $\theta_\epsilon$ is a subset of the parameter space identified in the proof. Also
\[ \lim_{\epsilon\rightarrow 0} {\rm d}_H \left( \theta_\epsilon, \theta_0 \right)= 0,\ \ {\rm in\ probability} \]
where ${\rm d}_H \left(\cdot,\cdot \right)$ is the Hausdorff distance.
\end{theorem}

\begin{theorem}[Drift estimation, homogenization problem]
\label{thm: drift hom}
Suppose that $\bar{f}_1$ in (\ref{homogenized eq}) depends on unknown parameters $\theta$, i.e. $\bar{f}_1(x) = \bar{f}_1(x;\theta)$. Let $\hat{\theta}(x;N,\delta)$ be the maximizer of the discretized likelihood corresponding to equation (\ref{averaged eq}) with step $\delta$, where $T = N\delta$. Suppose that we observe $\{ X_t, t\in[0,T] \}$ of system (\ref{eq:hom}) corresponding to $\theta = \theta_0$. Then, under appropriate assumptions described in \cite{PPS} (Theorem 4.5) and for $\delta = \epsilon^\alpha$ with $\alpha\in (0,2)$ and $N = [\epsilon^{-\gamma}]$ for $\gamma >\alpha$, it is possible to show that
\[ \lim_{\epsilon\rightarrow 0}  \hat{\theta}(X;N,\delta)= 0,\ \ {\rm in\ probability}. \]
\end{theorem}

The next two theorems deal with the estimation of the diffusion parameter of the limiting equation, given that this is constant. In that case, the MLE is the Quadratic Variation of the process. They assume that the dimension of the slow variable is $1$.

\begin{theorem}[Diffusion estimation, averaging problem]
\label{thm: diffusion av}
Let $X$ be the solution of (\ref{eq:av}) for $\bar{\sigma_1}\equiv \theta$ a constant. Then, under appropriate conditions described in \cite{PS} (Theorem 3.4) and for every $\epsilon>0$, we have that
\[\lim_{\delta\rightarrow 0}\frac{1}{N\delta}\sum_{n=0}^{N-1} | X_{(n+1)\delta} - X_{n\delta} |^2 = \theta^2\ \ {\rm a.s.}
\]
where $T = N\delta$ is fixed.
\end{theorem}

\begin{theorem}[Diffusion estimation, averaging problem]
\label{thm: diffusion hom}
Let $X$ be the solution of (\ref{eq:hom}) for $\bar{\tau}\equiv \theta$ a constant. Then, under appropriate conditions described in \cite{PS} (Theorem 3.5) and for  $\delta = \epsilon^\alpha$ with $\alpha\in(0,1)$, we have that
\[\lim_{\delta\rightarrow 0}\frac{1}{N\delta}\sum_{n=0}^{N-1} | X_{(n+1)\delta} - X_{n\delta} |^2 = \theta^2\ \ {\rm a.s.}
\]
where $T = N\delta$ is fixed.
\end{theorem}
It is conjecture that Theorem \ref{thm: diffusion hom} should hold for any $\alpha \in (0,2)$ and that the optimal $\alpha$, i.e. the one that minimizes the error, is $\alpha = \frac{2}{3}$.

Clearly, the most interesting case is that of estimating the diffusion parameter of the homogenized system. This is the case that we will study in detail in section \ref{my stuff}, assuming that the process is an Ornstein-Uhlenbeck process. Also, note that when estimating the diffusion parameter, the length of the time interval $T$ is fixed. We will relax this condition later on, for reasons explained in the following section.

\subsection{Equation-free modeling}

In practical applications it is often the case that the limiting equations (\ref{averaged eq}) and (\ref{homogenized eq}) are completely unknown. More generally, let us say that we have good reasons to believe that a certain variable of a multiscale system that evolves slowly behaves like a diffusion at a certain scale but we have complete ignorance of its drift and diffusion coefficients. We would like to find a way to estimate these coefficients. In statistical terms, let us say that we are interested in the non-parametric estimation of the drift and diffusion coefficients of the limiting equation. Note that our data comes ``on demand'' but for a certain cost, by simulating the multiscale model for given conditions.

A general algorithm for answering questions regarding the limiting dynamics of a quantity coming from a multiscale system that evolves slowly, when these are not explicitly known, comes under the name of ``equation-free'' algorithm (see \cite{Kevrekidis}). In our case, this would suggest pairing the problem of local estimation with an interpolation algorithm in order to estimate the drift and diffusion functions, denoted by $\bar{f}(x)$ and $\bar{\sigma}(x)$ respectively. We make this more concrete by describing the corresponding algorithm:
\begin{itemize}
\item[0.] Choose some initial condition $x_0$ and approximate $\bar{f}(x)$ and $\bar{\sigma}(x)$ by a local (polynomial) approximation around $x_0$. Simulate short paths of the multiscale system, so that the local approximation is acceptable. Note that the smaller the path, the better or simpler the local approximation.
\item[1.] For $n\geq 1$, choose another starting point $x_n$ using the knowledge of $\bar{f}(x_{n-1})$ and $\bar{\sigma}(x_{n-1})$ and possibly some of their derivatives on $x_{n-1}$, according to the rules of your interpolation algorithm. 
\item[2.] Repeat step 0, replacing $x_0$ by $x_n$.
\end{itemize}
As mentioned above, the size of the path $T$ needs to be small and possibly comparable to $\epsilon$. This is what led us to consider the estimation problem for $T=\epsilon^\alpha$. 

\section{The p-variation estimate}
\label{my stuff}

In this section, we study the problem of estimating the diffusion parameter of the homogenization limit of a simple multiscale Ornstein-Uhlenbeck process. We hope that the detailed analysis will provide some intuition for the general problem.

Consider the following system:
\begin{equation}
\label{system}
\begin{array}{ccc}
dY^{1,\epsilon}_t &=& \frac{\sigma}{\epsilon}Y^{2,\epsilon}_t dt \\
&&\\
dY^{2,\epsilon}_t &=& -\frac{1}{\epsilon^2}Y^{2,\epsilon}_t dt + \frac{1}{\epsilon}dW_t
\end{array}
\end{equation}
with initial conditions $Y^{1,\epsilon}_0 = y_1$ and $Y^{2,\epsilon}_0 = y_2$. It is not hard to see that the homogenization limit as $\epsilon\rightarrow 0$ is
\[ Y^{1,\epsilon}_t \rightarrow y_1 + \sigma W_t \]
and the convergence holds in a strong sense:
\begin{equation}
\label{L1 conv} 
\sup_{t\in[0,T]} |Y^{1,\epsilon}_t - y_1 - \sigma W_t| \stackrel{L_1}{\rightarrow} 0,\ \ {\rm as}\ \epsilon\rightarrow 0.
\end{equation}
Note that for this particular example, $Y^{1,\epsilon}_t$ is exactly equal to
\[ Y^{1,\epsilon}_t = y_1 + \sigma W_t - \epsilon\sigma\left( Y^{2,\epsilon}_t - y_2 \right) \]
and thus proving (\ref{L1 conv}) is equivalent to proving that
\[ \epsilon \sup_{t\in[0,T]} |Y^{2,\epsilon}_t - y_2| \stackrel{L_1}{\rightarrow} 0,\ \ {\rm as}\ \epsilon\rightarrow 0.\]
This follows from \cite{GP}.

We want to estimate the diffusion parameter $\sigma$ given a path $\{ Y^{1,\epsilon}_t(\omega)\ ;t\in [0,T] \}$. If we were to follow the approach discussed in the previous section, we would use the maximum likelihood estimate that corresponds to the limiting equation. In this case, this would be the quadratic variation. However, as discussed earlier, this is not a good estimate since the quadratic variation for any fixed $\epsilon>0$ is zero. To correct this, we subsample the data, which leads to the following estimate:
\begin{equation}
\label{QV estimate} 
\hat{\sigma}^2_\delta = \frac{1}{N\delta}\sum_{i=1}^N \left( Y^{1,\epsilon}_{i\delta} - Y^{1,\epsilon}_{(i-1)\delta} \right)^2,\ \ {\rm for}\ N=\frac{T}{\delta} 
\end{equation}
The asymptotic behavior of this estimate has been studied in \cite{PS, ABT}. In fact, taking advantage of the simplicity of the model, we can compute the $L_2$-error exactly, as a function of $\delta, \epsilon$ and $N$. We find that 
\begin{eqnarray}
\nonumber \frac{1}{\sigma^4}{\mathbb E}\left( \hat{\sigma}^2_\delta - \sigma^2 \right)^2 &=& 
\frac{\epsilon^4}{\delta^2}\left(1-e^{-\frac{\delta}{\epsilon^2}}\right)^2 \\
&+& 
\nonumber \left(2-4\frac{\epsilon^2}{\delta}\left(1-e^{-\frac{\delta}{\epsilon^2}}\right)+\frac{\epsilon^4}{\delta^2}\left(1-e^{-\frac{\delta}{\epsilon^2}}\right)^2\frac{3+e^{-\frac{\delta}{\epsilon^2}}}{1+e^{-\frac{\delta}{\epsilon^2}}}\right)\left( \frac{1}{N} \right) \\
&+&\frac{\epsilon^4}{\delta^2}\left(\frac{1-e^{-\frac{\delta}{\epsilon^2}}}{1+e^{-\frac{\delta}{\epsilon^2}}}\right)^2\left( \frac{e^{-\frac{2\delta N}{\epsilon^2}}-1}{N^2} \right)
\end{eqnarray}
For reasons explained earlier, we are interested in the behavior of this error not only when $T$ is fixed but also for $T\rightarrow 0$. Thus, we set $T=\epsilon^\alpha$ and, as before, $\delta = \epsilon^{\alpha+\beta}$, which lead to $N=\epsilon^{-\beta}$. We are interested in the behavior of the error as $\epsilon\rightarrow 0$. For these choices of $T$ and $\delta$, the square error will be
\begin{equation}
{\mathbb E}\left( \hat{\sigma}^2_\delta - \sigma^2 \right)^2 \sim {\mathcal O}\left( \epsilon^{4-2(\alpha+\beta)} + \epsilon^{2-\alpha} + \epsilon^{\beta}\right)
\end{equation}
For $\alpha$ fixed, we see that the error will be small if $0<\beta<2-\alpha$. In fact, the optimal choice for $\beta$ is $\beta = \frac{4-2\alpha}{3}$, in which case the error becomes
\begin{equation}
\left({\mathbb E}\left( \hat{\sigma}^2_\delta - \sigma^2 \right)^2 \right)^{\frac{1}{2}} \sim {\mathcal O}\left(  \epsilon^{\frac{2-\alpha}{3}} \right)
\end{equation}
So, for $\alpha=0$, we get that the optimal sub sampling rate is $\beta = \frac{4}{3}$, which results to an optimal error of order ${\mathcal O}\left(  \epsilon^{\frac{2}{3}} \right)$. However, if $\alpha>0$, the error can increase significantly, especially for non-optimal choices of $\delta$. 

In the rest of this section, we are going to investigate the behavior of the $p$-variation norm as an estimator of $\sigma$. The intuition comes from the following observation: we know that at scale ${\mathcal O}\left(  1 \right)$, $\{ Y^{1,\epsilon}_t(\omega)\ ;t\in [0,T] \}$ behaves like scaled Brownian motion while at scale ${\mathcal O}\left(  \epsilon \right)$, it is a process of bounded variation (finite length). Could it be that at scale ${\mathcal O}\left(  \epsilon^\alpha \right)$, the process behaves like a process of finite $p$-variation, for some $p$ that depends on $\alpha$? If so, would the $p$-variation norm be a better estimator of $\sigma$?

\subsection{The total $p$-variation}

We say that a real-valued continuous path $X:[0,T] \rightarrow {\mathbb R}$ has finite total $p$-variation if
\begin{equation}
\label{p-variation norm}
D_p\left(X \right)_T := \sup_{{\mathcal D}\left([0,T]\right)}\left(\sum_{t_\ell\in{\mathcal D}\left([0,T]\right)} |X_{t_{\ell+1}}-X_{t_\ell}|^p\right)^\frac{1}{p} <+\infty,
\end{equation}
where ${\mathcal D}\left([0,T]\right)$ goes through the set of all finite partitions of the interval $[0,T]$ (see also \cite{Terry book}). It is clear by the definition that a process of bounded variation will always have finite total $p$-variation for any $p>1$. Also, note that the total $p$-variation as defined above will only be zero if the process is constant. Thus, the total $p$-variation of a non-constant bounded variation process will always be a positive number. 

For $\epsilon>0$ fixed, the process $Y^{1,\epsilon}:[0,T] \rightarrow {\mathbb R}$ defined in (\ref{system}) is clearly of bounded variation, but its total variation is of order ${\mathcal O}\left( \frac{T}{\epsilon} \right)$. We will say that at scale ${\mathcal O}\left(  \epsilon^\alpha \right)$, the process $Y^{1,\epsilon}$ behaves like a process of finite total $p$-variation {\it in the limit} if
\begin{equation}
\lim_{\epsilon\rightarrow 0} \left( D_p\left(Y^{1,\epsilon} \right)_{\epsilon^\alpha}\right) <+\infty\ \ {\rm and}\ \ \forall q<p,\ \ \lim_{\epsilon\rightarrow 0} \left( D_q\left(Y^{1,\epsilon} \right)_{\epsilon^\alpha}\right) = +\infty.
\end{equation}
We will prove the following:
\begin{theorem}
At scale ${\mathcal O}\left(  \epsilon^\alpha \right)$ and $1<\alpha<2$, the process $Y^{1,\epsilon}:[0,T] \rightarrow {\mathbb R}$ defined in (\ref{system}) behaves like a process of finite total $(2-\alpha)$-variation in the limit.
\end{theorem}
First, we prove the following lemma:
\begin{lemma}
Let $X:[0,T] \rightarrow {\mathbb R}$ be a real-valued differentiable path of bounded variation. Then, its total $p$-variation is given by
\begin{equation}
\label{exact p-variation norm}
D_p\left(X \right)_T := \sup_{{\mathcal E}\left([0,T]\right)}\left(\sum_{t_\ell\in{\mathcal E}\left([0,T]\right)} |X_{t_{\ell+1}}-X_{t_\ell}|^p\right)^\frac{1}{p},
\end{equation}
where ${\mathcal E}\left([0,T]\right)$ goes through all finite sets of extremals of $X$ in the interval $[0,T]$. \end{lemma}
\begin{proof}
Consider the function
\[ f_{a,b}(t) = |X_t - X_a|^p + |X_b - X_t|^p,\ \ \ a<t<b. \]
This is maximized for $t$ an extremal point ($\dot{X}_t = 0$) or at $t = a$ or $t=b$. Thus, if ${\mathcal D} = \{0,t_1,\dots,t_{n-1},t_n = T\}$, there exists a set of extremals ${\mathcal E}$ with cardinality $|{\mathcal E}|\leq n+1$, such that
\[ \sum_{t_\ell\in{\mathcal D}} |X_{t_{\ell+1}}-X_{t_\ell}|^p \leq \sum_{t_\ell\in{\mathcal E}} |X_{t_{\ell+1}}-X_{t_\ell}|^p.\]
The set ${\mathcal E}$ can be constructed by choosing $\tau_1$ so that $f_{0,t_2}(t)$ is maximized and $\tau_k$ so that $f_{\tau_{k-1},t_{k+1}}(t)$ is maximized, for $k=2,\dots,n-1$. Thus
\[ \sup_{{\mathcal D}\left([0,T]\right)}\left(\sum_{t_\ell\in{\mathcal D}\left([0,T]\right)} |X_{t_{\ell+1}}-X_{t_\ell}|^p\right)^\frac{1}{p} \leq \sup_{{\mathcal E}\left([0,T]\right)}\left(\sum_{t_\ell\in{\mathcal E}\left([0,T]\right)} |X_{t_{\ell+1}}-X_{t_\ell}|^p\right)^\frac{1}{p}.\]
The opposite inequality is obvious and completes the proof.
\end{proof}

To prove the theorem, first we notice that 
\begin{equation}
\label{simple Dp} 
D_p(Y^{1,\epsilon})_T = \epsilon \sigma D_p(Z^{1})_{\frac{T}{\epsilon^2}},
\end{equation}
where $(Z^1,Z^2)$ satisfy
\begin{eqnarray*}
dZ^{1}_t &=& Z^{2}_t dt \\
dZ^{2}_t &=& -Z^{2}_t dt + dW_t
\end{eqnarray*}
Now, $Z^1$ is clearly differentiable and thus, by the lemma 
\[ D_p(Z^1)_T = \sup_{{\mathcal E}\left([0,T]\right)}\left(\sum_{t_\ell\in{\mathcal E}\left([0,T]\right)} |Z^1_{t_{\ell+1}}-Z^1_{t_\ell}|^p\right)^\frac{1}{p}
\]
The derivative of $Z^1$ is equal to $Z^2$, so all its extremal points correspond to zero-crossings of $Z^2$. So, for $s,t\in{\mathcal E}$, 
\[ Z^1_t - Z^1_s = \left( W_t-W_s \right) - \left( Z^2_t - Z^2_s \right) = W_t - W_s\]
and $D_p(Z^1)_T$ becomes
\begin{eqnarray}
\label{DpZ}
\nonumber D_p(Z^1)_T &=& \sup_{{\mathcal E}\left([0,T]\right)}\left(\sum_{t_\ell\in{\mathcal E}\left([0,T]\right)} |W_{t_{\ell+1}}-W_{t_\ell}|^p\right)^\frac{1}{p} = \\
&=& \lim_{\delta\rightarrow 0}\left(\sum_{t_\ell\in{\mathcal E_\delta}\left([0,T]\right)} |W_{t_{\ell+1}}-W_{t_\ell}|^p\right)^\frac{1}{p},
\end{eqnarray}
where 
\[{\mathcal E_\delta}\left([0,T]\right) = \{ 0 = t_0,t_1,\dots, t_{N_\delta(T)},T \}\]
and $\{ t_1,\dots, t_{N_\delta(T)}\}$ is the set of all zero-crossings of $Z^2$ in $[0,T]$ that are at least distance $\delta$ apart from each other, i.e. if $t_k\in {\mathcal E_\delta}\left([0,T]\right)$ and $k<N_\delta(T)$, then $t_{k+1}$ is the first time that $Z^2$ crosses zero after time $t_k +\delta$. Note that the set of zero-crossings of $Z^2$ in $[0,T]$ is an uncountable set that contains no intervals with probability 1. Equation (\ref{DpZ}) follows from the following two facts: (i) in general, adding any point to the partition will increase the $L_p$ norm and thus the supremum is achieved for a countable set of zero-crossings and (ii) any countable set that is dense in the set of all zero-crossings will give the same result.

If $\tau_\delta$ is the stopping time of the first zero crossing of $Z^2$ after $\delta$ given $Z^2_0 = 0$, then the random variables $\{ \tau^\delta_k = (t_{k}-t_{k-1}), t_k \in {\mathcal E}_\delta\left( [0,T]\right), k\leq N_\delta(T)\}$ are i.i.d. with the same law as that of $\tau_\delta$. Thus, the sum $\sum_{t_\ell\in{\mathcal E_\delta}\left([0,T]\right)} |W_{t_{\ell+1}}-W_{t_\ell}|^p$ is a sum of i.i.d. random variables of finite mean (to be computed in the following section) and as a consequence of the Law of Large Numbers, it grows like $N_\delta(T)$. From \cite{F-Z}, we know that $N_\delta(T)\sim{\mathcal O}\left( \frac{T}{\epsilon} \right)$. We conclude that
\[ D_p(Z^1)_T \sim {\mathcal O}\left( T^\frac{1}{p} \right).\]
Finally, from (\ref{simple Dp}), it is clear that
\[ D_p(Y^{1,\epsilon})_{\epsilon^\alpha} \sim {\mathcal O}\left( \epsilon\left(\frac{\epsilon^\alpha}{\epsilon^2}\right)^\frac{1}{p} \right) \sim {\mathcal O}\left( \epsilon^{1+\frac{\alpha-2}{p}} \right),\]
which proves the theorem.

\subsection{The $p$-variation estimates}

Similar to the quadratic variation estimate $\hat{\sigma}^2$ defined in (\ref{QV estimate}), we define the $p$-variation estimates as the properly normalized total $p$-variation of the process:
\begin{equation}
\label{p-var estimate}
\hat{\sigma}^p := \frac{1}{C_p(T)}\left(D_p(Y^{1,\epsilon})_T\right)^p.
\end{equation}
We will study the $L_2$-error of this estimate in different scales. First, we need to define the constant $C_p(T)$. The natural choice would be to choose $C_p(T)$ so that ${\mathbb E}\left( \hat{\sigma}^p\right) = \sigma^p$. So,
\[ C_p(T) = \frac{1}{\sigma^p} {\mathbb E}\left( \left( D_p(Y^{1,\epsilon})_T\right)^p \right). \]
We need to compute ${\mathbb E}\left( \left( D_p(Y^{1,\epsilon})_T\right)^p \right)$. From (\ref{simple Dp}), we get that
\[ {\mathbb E}\left( \left( D_p(Y^{1,\epsilon})_T\right)^p \right) = \epsilon^p \sigma^p {\mathbb E}\left( \left( D_p(Z^1)_{\frac{T}{\epsilon^2}}\right)^p \right).\]
Using (\ref{DpZ}), we get that
\[
{\mathbb E}\left( \left( D_p(Z^1)_{T}\right)^p \right) = \lim_{\delta\rightarrow 0}{\mathbb E} \left(\sum_{t_\ell\in{\mathcal E_\delta}\left([0,T]\right)} |W_{t_{\ell+1}}-W_{t_\ell}|^p\right), 
\]
Note that for any $p>1$, $D_p(Z^1)_{T}\leq D_1(Z^1)_{T}$, where ${\mathbb E}\left( D_1(Z^1)_{T}^p\right) < +\infty$. Thus, from the Dominated Convergence Theorem, the limit can come out of the expectation. To simplify our computations, from now on we will assume that $Z^2_0 = Z^2_T = 0$. We have already observed that the random variables $\left\{(W_{t_{\ell+1}}-W_{t_\ell}), t_\ell\in{\mathcal E_\delta}\left([0,T]\right), \ell< N_\delta(T) \right\}$ are independent and distributed like $W_{\tau_\delta}$ where $\tau_\delta$ is the first time $Z^2$ crosses zero after $t=\delta$, given that $Z^2_0 = 0$.  Thus,
\[ {\mathbb E} \left(\sum_{t_\ell\in{\mathcal E_\delta}\left([0,T]\right)} |W_{t_{\ell+1}}-W_{t_\ell}|^p\right) = {\mathbb E} N_\delta(T)\ {\mathbb E} |W_{\tau_\delta}|^p + {\mathbb E}|W_T - W_{t_{N_\delta(T)}}|^p,
\]
where $N_\delta(T)$ is the number of zero-crossings of $Z^2$ in interval $[0,T]$ that are distance $\delta$ apart from each other. First, we notice that 
\[ {\mathbb E} |W_{\tau_\delta}|^p = {\mathbb E} \left({\mathbb E}\left( |W_{\tau_\delta}|^p  {\big |} \tau_\delta \right)\right) = \frac{1}{\sqrt{\pi}}2^\frac{p}{2}\Gamma\left(\frac{p+1}{2}\right) {\mathbb E} \left((\tau_\delta)^\frac{p}{2}\right).
\]
To compute ${\mathbb E} \left((\tau_\delta)^p\right)$, we note that $\tau_\delta$ can be written as $\tau_\delta = \delta + \tau(Z^2_\delta)$, where $\tau(z)$ is the first zero-crossing of the process $Z^2$ given that it starts at $z$. For $Z^2$ an Ornstein-Uhlenbeck process, the p.d.f. of $\tau(z)$ has been computed explicitly (see \cite{Ricciardi-Sato}) and is given by 
\[ f(t,z) = \frac{2}{\sqrt{\pi}}\frac{|z| e^{-t}}{(1-e^{-2t})^\frac{3}{2}}\exp\left(-\frac{z^2 e^{-2t}}{1-e^{-2t}}\right).
\]
Since $Z^2_0 = 0$ by assumption, $Z^2_\delta$ is a Gaussian random variable with zero mean and variance $\frac{1}{2}(1-e^{-2\delta})$. Let us denote its p.d.f. by $g_\delta(z)$. It follows that the p.d.f. of $\tau(Z^2_\delta)$ is given by
\begin{equation}
\label{pdf tau} 
h_\delta(t) := \int_{-\infty}^\infty f(t,z)g_\delta(z)dz = \frac{4 e^t {\rm csch}(\delta + t)\sinh(\delta)\sinh(t)}{\sqrt{(1-e^{-2\delta})(1-e^{-2t})}(-1+e^{2t})\pi}
\end{equation}
where ${\rm csch}(t)=\frac{1}{\sinh(t)}$ and $\sinh(t)$ the hyperbolic sine. We write 
\begin{eqnarray*}
\lim_{\delta\rightarrow 0} \frac{1}{\sqrt{\delta}} {\mathbb E} \left((\tau_\delta)^p\right) &=& \lim_{\delta\rightarrow 0} \frac{1}{\sqrt{\delta}}\int_0^\infty (\delta+t)^p h_\delta(t) dt = \\
&=& \lim_{\delta\rightarrow 0} \int_0^\infty \frac{(\delta+t)^p h_\delta(t)}{\sqrt{\delta} t H(t)}  t H(t) dt.
\end{eqnarray*}
where 
\[ H(t) = \frac{4e^{-t}\sqrt{e^{-t} \sinh(t)}}{(1-e^{-2t})^2 \pi}\ \ {\rm and}\ \ \int_0^\infty t H(t) dt = \sqrt{2}.
\]
The function $\frac{(\delta+t)^p h_\delta(t)}{\sqrt{\delta} t H(t)}$ is increasing to $t^{-1+p}$ as $\delta \downarrow 0$ and thus by the dominated theorem we find that 
\begin{equation}
\label{Kp}
K_p := \int_0^\infty t^p H(t) dt
\end{equation}
Notice that for $t\rightarrow 0$, $H(t)$ behaves like $t^{-\frac{3}{2}}$ and thus the integral $K_p$ is finite if and only if $p>\frac{1}{2}$. Also, for $p=1,2$ we find that $K_1 = \sqrt{2}$ and $K_2 = 2\sqrt{2}\log{2}$.


Now, we need to compute the limit of $\sqrt{\delta} {\mathbb E} N_\delta(T)$ as $\delta\rightarrow 0$. We can use the results in \cite{F-Z} to get an upper and lower bound and show that $N_\delta(T)$  behaves like ${\mathcal O}\left( \frac{T}{\sqrt{\delta}}\right)$. However, we need to know the exact value of the limit. We proceed as follows: we write
\begin{equation}
\label{EN 1}
{\mathbb E} N_\delta(T) = \sum_{n=1}^\infty {\mathbb P}\left( N_\delta(T)\geq n \right) = \sum_{n=1}^\infty {\mathbb P}\left( \sum_{i=1}^n \tau_i^\delta \leq T \right) 
\end{equation}
where $\tau_i^\delta = t_i-t_{i-1}$ for $t_i \in {\mathcal E}_\delta \left( [0,T] \right)$ and $i\leq N_\delta(T)$. Using (\ref{pdf tau}) we find that the Laplace transform of the distribution of $\tau_\delta$ is 
\begin{equation}
\label{Laplace 1}
\hat{H}_\delta \left( \lambda \right) = e^{- \lambda \delta} \hat{h}_\delta \left( \lambda \right) = \frac{2 e^{-d(\lambda +1)} \Gamma(\frac{\lambda + 1}{2}) \sinh(d)}{\sqrt{\pi(1-e^{-2d})}} \bar{F}_1\left( 1, \frac{\lambda + 1}{2},\frac{\lambda + 2}{2},e^{-2 d}\right),
\end{equation}
where $\bar{F}_1\left( a,b,c,x\right)$ is the regularized hypergeometric function given by
\[ \bar{F}_1\left( a,b,c,x\right) = \frac{1}{\Gamma(c)}\sum_{k=0}^\infty \frac{(a)_k (b)_k}{(c)_k} \frac{z^k}{k!},\ \ {\rm and}\ \ (d)_n = \prod_{k=0}^{n-1}(d+k).
\]
We find that for small $d>0$, this behaves like 
\begin{equation}
\label{Laplace 1 expansion}
\hat{H}_\delta \left( \lambda \right) = 1 -  \frac{2\sqrt{2}}{\sqrt{\pi}}\frac{\Gamma(\frac{\lambda+1}{2})}{\Gamma(\frac{\lambda}{2})}\sqrt{d} + {\mathcal O}\left(\delta\right).
\end{equation}
Since the $\tau_i^\delta$'s are i.i.d., the Laplace transform of the sum $\sum_{i=1}^n \tau_i^\delta$ will be $ \hat{H}_\delta \left( \lambda \right)^n$ and thus we write
\[ {\mathbb P}\left( \sum_{i=1}^n \tau_i^\delta \leq T \right)  = \int_0^T {\mathcal L}^{-1}[\hat{H}_\delta \left( \lambda \right)^n](dt),
\]
where ${\mathcal L}^{-1}$ denotes the operator of the inverse Laplace transform. Substituting this back to (\ref{EN 1}), we get
\begin{eqnarray}
\label{EN 2}
\nonumber {\mathbb E} N_\delta(T) &=& 
\sum_{n=1}^\infty \int_0^T {\mathcal L}^{-1}[\hat{H}_\delta \left( \lambda \right)^n](dt) \\ 
\nonumber &=& \int_0^T {\mathcal L}^{-1}[\sum_{n=1}^\infty \hat{H}_\delta \left( \lambda \right)^n](dt) \\
&=& \int_0^T {\mathcal L}^{-1}[\frac{\hat{H}_\delta \left( \lambda \right)}{1-\hat{H}_\delta \left( \lambda \right)}](dt).
\end{eqnarray}
Taking the limit inside the operator, we finally see that
\begin{equation}
\label{lim EN}
\lim_{\delta\rightarrow\infty} \sqrt{\delta} {\mathbb E}\left( N_\delta(T) \right) = \frac{\sqrt{\pi}}{2\sqrt{2}}\int_0^T {\mathcal L}^{-1}[\frac{\Gamma(\frac{\lambda}{2})}{\Gamma(\frac{\lambda+1}{2})}](dt) = \frac{T}{\sqrt{2}}.
\end{equation}
Finally, we note that since $Z^2_T = 0$ by assumption,
\[ (T-t_{N_\delta(T)}) < \delta \Rightarrow \lim_{\delta\rightarrow 0} {\mathbb E}|W_T - W_{t_{N_\delta(T)}}|^p = 0.
\]
For every $p>1$, we set 
\begin{equation}
\label{constants}
a_p := \frac{1}{\sqrt{\pi}}2^\frac{p}{2}\Gamma\left(\frac{p+1}{2}\right)\ \ {\rm and}\ \ c_p := \frac{a_p}{\sqrt{2}} K_{\frac{p}{2}}.
\end{equation}
Putting everything together, we find that
\begin{equation}
\label{moment 1}
{\mathbb E}\left( \left( D_p(Z^{1})_T\right)^p \right) := c_p T 
\end{equation}
and consequently
\[ {\mathbb E}\left( \left( D_p(Y^{1,\epsilon})_T\right)^p \right) = \epsilon^p \sigma^p c_p \frac{T}{\epsilon^2} = \epsilon^{p-2} \sigma^p c_p T.\]
Thus we set 
\begin{equation}
\label{Cp}
C_p(T) := \epsilon^{p-2} c_p T.
\end{equation}
By construction, the $p$-variation estimates $\hat{\sigma}^p$ defined in (\ref{p-var estimate}) are consistent, i.e. ${\mathbb E}\left( \hat{\sigma}^p\right) = \sigma^p$. We now compute its square $L_2$-error:
\begin{eqnarray}
\label{p-var error}
\nonumber {\mathbb E}\left( \hat{\sigma}^p - \sigma^p \right)^2 &=& {\mathbb E}\left( \frac{\left(D_p(Y^{1,\epsilon})_T\right)^p}{C_p(T)} - \sigma^p \right)^2 = \\
\nonumber &=& {\mathbb E}\left( \frac{\left(D_p(Y^{1,\epsilon})_T\right)^{2p}}{C_p(T)^2} \right) - \sigma^{2p} = \\
\nonumber&=& \frac{1}{C_p(T)^2}{\mathbb E}\left( \left(D_p(Y^{1,\epsilon})_T\right)^{2p} \right) - \sigma^{2p} = \\
\nonumber &=& \frac{\epsilon^{2p}\sigma^{2p}}{\epsilon^{2p-4} c_p^2 T^2}{\mathbb E}\left( \left(D_p(Z^1)_{\frac{T}{\epsilon^2}}\right)^{2p} \right) - \sigma^{2p} \\
&=& \sigma^{2p}\left( \frac{\epsilon^4}{c_p^2 T^2}{\mathbb E}\left( \left(D_p(Z^1)_{\frac{T}{\epsilon^2}}\right)^{2p}\right) - 1 \right)
\end{eqnarray}
To proceed, we need to compute the second moment of $\left( D_p(Z^1)_T \right)^p$. As with the computation of the first moment, we write:
\begin{eqnarray*}
{\mathbb E}\left( \left( D_p(Z^1)_{T}\right)^{2p} \right) &=& \lim_{\delta\rightarrow 0}{\mathbb E} \left(\sum_{t_\ell\in{\mathcal E_\delta}\left([0,T]\right)} |W_{t_{\ell+1}}-W_{t_\ell}|^p\right)^2 \\
&=& \lim_{\delta\rightarrow 0}{\mathbb E} \left(\sum_{n=1}^{N_\delta(T)} |W_{\tau^\delta_n}|^p + |W_{T}-W_{t_{N_\delta(T)}}|^p\right)^2 \\
&=& \lim_{\delta\rightarrow 0}{\mathbb E} \left(\sum_{n=1}^{N_\delta(T)} |W_{\tau^\delta_n}|^p\right)^2,
\end{eqnarray*}
where the last line comes from the fact that $\left( T-t_{N_\delta(T)} \right) < \delta$. To compute the above expectation, we write
\begin{eqnarray*}
{\mathbb E} \left(\sum_{n=1}^{N_\delta(T)} |W_{\tau^\delta_n}|^p\right)^2 &=& {\mathbb E}\left(\sum_{m,n=1}^{N_\delta(T)} |W_{\tau^\delta_m}|^p|W_{\tau^\delta_n}|^p\right)^2 \\
&=& {\mathbb E} N_\delta(T) {\mathbb E}|W_{\tau_\delta}|^{2p} + {\mathbb E} \left( N_\delta(T)^2 - N_\delta(T) \right)\left({\mathbb E}|W_{\tau_\delta}|^{p}\right)^2 \\
&=& {\mathbb E} N_\delta(T) {\mathbb E}|W_{\tau_\delta}|^{2p} + {\mathbb E} N_\delta(T)^2 \left({\mathbb E}|W_{\tau_\delta}|^{p}\right)^2 + {\mathcal O}\left(\sqrt{d}\right)
\end{eqnarray*}
where the last line follows from the fact that $N_\delta(T)\sim{\mathcal O}\left(\frac{T}{\sqrt{d}}\right)$ and ${\mathbb E}|W_{\tau_\delta}|^{p}\sim{\mathcal O}\left(\sqrt{\delta}\right)$. It remains to compute the limit of $\delta {\mathbb E}N_\delta(T)^2$. Following a similar approach to the one before, we write
\begin{eqnarray*}
{\mathbb E} N_\delta(T)^2 = \sum_{n=1}^\infty (2n-1){\mathbb P}\left( N_\delta(T)\geq n \right) = \sum_{n=1}^\infty (2n-1) {\mathbb P}\left( \sum_{i=1}^n \tau_i^\delta \leq T \right) \\
= 2 \sum_{n=1}^\infty n {\mathbb P}\left( \sum_{i=1}^n \tau_i^\delta \leq T \right) + {\mathcal O}\left( \frac{1}{\sqrt{\delta}} \right)
\end{eqnarray*}
and
\begin{eqnarray*}
\sum_{n=1}^\infty n {\mathbb P}\left( \sum_{i=1}^n \tau_i^\delta \leq T \right) &=& 
\sum_{n=1}^\infty n \int_0^T {\mathcal L}^{-1}[\hat{H}_\delta \left( \lambda \right)^n](dt) \\  
&=& \int_0^T {\mathcal L}^{-1}[\sum_{n=1}^\infty n \hat{H}_\delta \left( \lambda \right)^n](dt) \\
&=& \int_0^T {\mathcal L}^{-1}[\frac{\hat{H}_\delta \left( \lambda \right)}{\left( 1-\hat{H}_\delta \left( \lambda \right)\right)^2}](dt).
\end{eqnarray*}
Taking the limit as $\delta \rightarrow 0$, we get
\begin{eqnarray*}
\lim_{\delta \rightarrow 0}\delta {\mathbb E} N_\delta(T)^2 &=& \lim_{\delta \rightarrow 0}2\int_0^T {\mathcal L}^{-1}[\frac{\delta\hat{H}_\delta \left( \lambda \right)}{\left( 1-\hat{H}_\delta \left( \lambda \right)\right)^2}](dt) \\
&=& \frac{\pi}{4}\int_0^T {\mathcal L}^{-1}[\left(\frac{\Gamma(\frac{\lambda}{2})}{\Gamma(\frac{\lambda+1}{2})}\right)^2](dt) \\
&=& \frac{T^2}{2} + \left(2 \log{2}\right) T
\end{eqnarray*}

Putting everything together, we get
\begin{eqnarray*}
{\mathbb E}\left( \left( D_p(Z^1)_{T}\right)^{2p} \right) &=& \lim_{\delta\rightarrow 0}\left(\sqrt{\delta}{\mathbb E} N_\delta(T) {\mathbb E}\frac{|W_{\tau_\delta}|^{2p}}{\sqrt{\delta}} + \delta{\mathbb E} N_\delta(T)^2 \left({\mathbb E}\frac{|W_{\tau_\delta}|^{p}}{\sqrt{\delta}}\right)^2\right) \\
&=& \frac{T}{\sqrt{2}} a_{2p}K_p +  \left(\frac{T^2}{2} + \left(2 \log{2}\right) T\right) \left(a_p K_{\frac{p}{2}}\right)^2 \\
&=& T c_{2p} + \left(T^2 + \left(4 \log{2}\right) T\right)\left(c_p \right)^2\\
&=& T^2\left(c_p \right)^2 + T\left(c_{2p} + \left(4 \log{2}\right)\left(c_p \right)^2 \right),
\end{eqnarray*}
where $a_p$ and $c_p$ are defined in (\ref{constants}) and $K_p$ is defined in (\ref{Kp}). Finally, we get
\begin{eqnarray}
\label{p-var error 2}
\nonumber {\mathbb E}\left( \hat{\sigma}^p - \sigma^p \right)^2 &=& \sigma^{2p}\left( \frac{\epsilon^4}{T^2 c_p^2}\left(\frac{T^2}{\epsilon^4}\left(c_p \right)^2 + \frac{T}{\epsilon^2}\left(c_{2p} + \left(4 \log{2}\right)\left(c_p \right)^2 \right)\right) - 1 \right) \\
&=& \sigma^{2p}\frac{\epsilon^2}{T} \left( \frac{c_{2p}}{c_p^2 } +  4\log{(2)}\right) = \sigma^{2p}\frac{\epsilon^2}{T} E(p),
\end{eqnarray}
where $E(p) = \frac{c_{2p}}{c_p^2 } +  4\log{(2)}$. This is an increasing function for $p\in [1,2]$ and 
\[ 4\log{2} =: E(1) \leq E(p) \leq E(2) := 10\log{2}, \ \ \forall p\in [1,2].\]
We summarize our conclusions in the following 
\begin{theorem}
The $L_2$-error of the estimator $\hat{\sigma}^p$ defined in (\ref{p-var estimate}) is described by (\ref{p-var error 2}). At scale ${\mathcal O}\left( \epsilon^\alpha \right)$, the error is of order ${\mathcal O}\left( \epsilon^{\frac{2-\alpha}{2}} \right)$.
\end{theorem}
We see that the performance of the estimators $\hat{\sigma}^p$ is the same for all $p>1$ and they outperform the $\hat{\sigma}^2_\delta$ estimator defined in (\ref{QV estimate}). In terms of the constant $E(p)$, the smaller the $p$, the smaller the error. However, there is a problem: except for scale ${\mathcal O}(1)$ ($\alpha = 0$), the normalizing constant $C_p$ depends on $\epsilon$, which will in general be unknown. We go on to define a new estimator that does not assume knowledge of $\epsilon$. 

\subsection{Estimating the scale separation variable $\epsilon$}

Suppose that $T < 1$ and $T = \epsilon^\alpha$ for some $\alpha>0$. We define the new estimator $\tilde{\sigma}^p$ similar to $\hat{\sigma}^p$, only use $c_p$ rather than $C_p$ as our normalization constant. Thus, we define
\begin{equation}
\label{new p-var estimate}
\tilde{\sigma}^p = \frac{1}{c_p}\left(D_p(Y^{1,\epsilon})_T\right)^p,
\end{equation}
where $c_p$ is defined in (\ref{constants}). Then
\begin{eqnarray}
\label{new p-var error}
\nonumber {\mathbb E}\left( \tilde{\sigma}^p - \sigma^p \right)^2 &=& {\mathbb E}\left( \frac{\left(D_p(Y^{1,\epsilon})_T\right)^p}{c_p} - \sigma^p \right)^2 = \\
\nonumber &=& \frac{\epsilon^{2p}\sigma^{2p}}{c_p^2}{\mathbb E}\left( \left(D_p(Z^1)_{\frac{T}{\epsilon^2}}\right)^{2p} \right) - 2\sigma^{p}\frac{\epsilon^{p}\sigma^{p}}{c_p}{\mathbb E}\left( \left(D_p(Z^1)_{\frac{T}{\epsilon^2}}\right)^{p} \right) + \sigma^{2p} \\
\nonumber &=& \frac{\epsilon^{2p}\sigma^{2p}}{c_p^2}\left(\frac{T^2}{\epsilon^4}\left(c_p \right)^2 + \frac{T}{\epsilon^2}\left(c_{2p} + \left(4 \log{2}\right)\left(c_p \right)^2 \right)\right) - 2\sigma^{p}\frac{\epsilon^{p}\sigma^{p}}{c_p}\left( c_p\frac{T}{\epsilon^2} \right) + \sigma^{2p} \\
&=& \sigma^{2p}\left(\frac{T^2}{\epsilon^{4-2p}} + \frac{T}{\epsilon^{2-2p}}\left(\frac{c_{2p}}{\left(c_p \right)^2} + \left(4 \log{2}\right) \right)- 2\left( \frac{T}{\epsilon^{2-p}} \right) + 1 \right)
\end{eqnarray}
and by substituting $T$ by $\epsilon^\alpha$ this becomes
\begin{equation}
\label{new p-var error 2}
{\mathbb E}\left( \hat{\sigma}^p - \sigma^p \right)^2 = \sigma^{2p}\left( \epsilon^{2p+2a-4}  + \epsilon^{2p+a-2} E(p) - 2\epsilon^{p+a-2}+ 1 \right).
\end{equation}
Thus, we get the following behavior:
\begin{itemize}
\item[(i)] For $p>2-\alpha$, the error is of order ${\mathcal O}\left(1 \right)$.
\item[(ii)] For $p=2-\alpha$, the error is well-behaved and of order ${\mathcal O}\left(\epsilon^\frac{2-\alpha}{2} \right)$.
\item[(iii)] For $p<2-\alpha$ and $\alpha<2$, the error explodes like ${\mathcal O}\left(\epsilon^{2p+2a-4} \right)$.
\end{itemize}
We conclude that the optimal estimator is $\hat{\sigma}^2$, since it does not assume knowledge of $\epsilon$ and the estimators $\tilde{\sigma}^p$ do not outperform it even for $p=2-\alpha$ (except that the constant $E(p)$ is smaller). However, the estimators $\tilde{\sigma}^p$ can be used to estimate the scale separation variable $\epsilon$. We set
\[ \hat{p} := \arg\min_{1<p<2} | \left(\tilde{\sigma}^p\right)^\frac{1}{p} - \left(\hat{\sigma}^2\right)^\frac{1}{2} |\]
and
\[ \hat{\alpha} := 2-\hat{p}. \]
Then, we estimate $\epsilon$ by
\[ \hat{\epsilon} := T^\frac{1}{\hat{\alpha}}.\]


\begin{thebibliography}{99}

\bibitem{ABT}
R. Azencott, A. Beri and I. Timofeyev. Adaptive subsampling for parametric estimation of Gaussian Diffusions. {\it Preprint}.

\bibitem{Bishwal}
J.P.N. Bishwal. {\it Parameter Estimation in Stochastic Differential Equations}. Lecture Notes in Mathematics vol. 1923, Springer, Berlin, 2008.

\bibitem{Chen book}
X. Chen. {\it Limit Theorems for Functionals of Ergodic Markov Chains with General State Space}. Memoirs of the AMS, vol. 129, No. 664, 1999.

\bibitem{F-Z}
D. Florens-Zmirou. Statistics on crossings of discretized diffusions and local time. {\it Stochastic Process. Appl.} 39: pp. 139--151, 1991.

\bibitem{FW book}
M.I. Freidlin, A.D. Wentzell. {\it Random perturbations of dynamical systems}. Springer, New York, 1998.

\bibitem{GP}
S.E. Graversen and G. Peskir. Maximal Inequalities for the Ornstein-Uhlenbeck Process. {\it Proceedings of the A.M.S.} 128(10): pp. 3035--3041, 2000.

\bibitem{KMS 1}
M. Katsoulakis, A. Majda and A. Sopasakis. Multiscale couplings in prototype hybrid deterministic/stochastic systems: Part 1, deterministic closures. {\it Comm. Math. Sci.} 2, pp. 255–-294, 2004

\bibitem{KMS 2}
M. Katsoulakis, A. Majda and A. Sopasakis. Multiscale couplings in prototype hybrid deterministic/stochastic systems: Part 2, stochastic closures. {\it Comm. Math. Sci.} 3, pp. 453–-478, 2005.

\bibitem{Kutoyants}
Y.A. Kutoyants. {\it Statistical inference for ergodic diffusion processes}. Springer-Verlag, London, 2004.

\bibitem{Kevrekidis}
J. Li, P.G. Kevrekidis, C.W. Gear and I.G. Kevrekidis. Deciding the nature of the coarse equation through microscopic simulations: the baby-bathwater scheme. {\it SIAM Review}, 49(3): pp. 469--487, 2007.

\bibitem{Terry book}
T. Lyons and Z. Qian. {\it System control and rough paths.} Oxford University Press, Oxford, 2002.

\bibitem{MTV 1}
A.J. Majda, I. Timofeyev and E. Vanden-Eijnden. A mathematics framework for stochastic climate models. {\it Comm. Pure Appl. Math.} 54, pp. 891–-974, 2001.

\bibitem{MTV 2}
A.J. Majda, I. Timofeyev and E. Vanden-Eijnden. Stochastic models for selected slow variables in large deterministic systems. {\it Nonlinearity} 19, pp. 769–-794, 2006.

\bibitem{PPS}
A. Papavasiliou, G.A. Pavliotis and A.M. Stuart. Maximum likelihood drift estimation for multiscale diffusions.  {\it Stochastic Process. Appl.}  119, pp. 3173--3210, 2009.

\bibitem{PS}
G.A. Pavliotis and A.M. Stuart. Parameter estimation for multiscale diffusions. {\it J. Stat. Phys.} 127, pp. 741--781, 2007.

\bibitem{PS book}
G.A. Pavliotis and A.M. Stuart. {\it Multiscale Methods: Averaging and Homogenization}. Texts in Applied Mathematics 53, Springer, New York, 2008.

\bibitem{Ricciardi-Sato}
L.M. Ricciardi and S. Sato. First-passage-time density and moments of the Ornstein-Uhlenbeck process. {\it J. Appl. Prob.} 25, pp. 43--57, 1988.

\end{thebibliography}
\end{document}